\documentclass[fleqn, a4paper, oneside]{imsart}

\usepackage{amsthm,amsmath,amssymb}
\usepackage{graphicx}
\usepackage{makeidx}
\usepackage[authoryear]{natbib}

\bibpunct{(}{)}{;}{a}{,}{,}
\RequirePackage[OT1]{fontenc}

\startlocaldefs

\newtheorem{Definition}{Definition}
\newtheorem{Proposition}{Proposition}
\newtheorem{Theorem}{Theorem}

\newtheorem{Lemma}{Lemma}
\newenvironment{Proof}{\noindent{\bf Proof:}}{\hfill\rule{2mm}{2mm}\vspace{0.3cm}}
\newenvironment{Remark}{\vspace{0.1cm}\noindent{\bf Remark.}}{\vspace{0.3cm}}

\theoremstyle{remark}

\theoremstyle{definition}

\def\R{\mathbb{R}}
\def\C{\mathbb{C}}
\def\N{\mathbb{N}}

\def\eps{\varepsilon}

\def\GG{\mathcal{G}}
\def\TT{\mathcal{T}}

\endlocaldefs

\begin{document}

\begin{frontmatter}


\title{Robust nonparametric detection of objects in noisy images.}
\runtitle{Robust detection in noisy images.}


\begin{aug}
\author{\snm{Mikhail} \fnms{Langovoy}\corref{}\thanksref{t2}
\ead[label=e1]{langovoy@eurandom.tue.nl}}

\affiliation{
        EURANDOM,\\
         The Netherlands.}

\address{Mikhail Langovoy, Technische Universiteit Eindhoven, \\
EURANDOM, P.O. Box 513,
\\
5600 MB, Eindhoven, The Netherlands\\
\printead{e1}\\
Phone: (+31) (40) 247 - 8113\\
Fax: (+31) (40) 247 - 8190\\}

\and

\author{\snm{Olaf} \fnms{Wittich} \ead[label=e2]{o.wittich@tue.nl}}

\affiliation{
        Technische Universiteit Eindhoven and EURANDOM,\\
         The Netherlands.}

\address{Olaf Wittich, Technische Universiteit Eindhoven and \\
EURANDOM, P.O. Box 513,
\\
5600 MB, Eindhoven, The Netherlands\\
\printead{e2}\\
Phone: (+31) (40) 247 - 2499}

\thankstext{t2}{Corresponding author.}

\runauthor{M. Langovoy and O. Wittich}
\end{aug}

\begin{abstract}
We propose a novel statistical hypothesis testing method for detection of objects in noisy images. The method uses results from percolation theory and random graph theory. We present an algorithm that allows to detect objects of unknown shapes in the presence of nonparametric noise of unknown level and of unknown distribution. No boundary shape constraints are imposed on the object, only a weak bulk condition for the object's interior is required. The algorithm has linear complexity and exponential accuracy and is appropriate for real-time systems. 

In this paper, we develop further the mathematical formalism of our method and explore important connections to the mathematical theory of percolation and statistical physics. We prove results on consistency and algorithmic complexity of our testing procedure. In addition, we address not only an asymptotic behavior of the method, but also a finite sample performance of our test.\\
\end{abstract}



\end{frontmatter}

\section{Introduction}\label{Section1}

Assume we observe a noisy digital image on a screen of $N \times N$ pixels. Object detection and image reconstruction for noisy images are two of the cornerstone problems in image analysis. In this paper, we propose a new efficient technique for quick detection of objects in noisy images. Our approach uses mathematical percolation theory.

Detection of objects in noisy images is the most basic problem of image analysis. Indeed, when one looks at a noisy image, the first question to ask is whether there is any object at all. This is also a primary question of interest in such diverse fields as, for example, cancer detection (\cite{Cancer_Detection_1}), automated urban analysis (\cite{Road_Detection_IEEE}), detection of cracks in buried pipes (\cite{Sinha200658}), and other possible applications in astronomy, electron microscopy and neurology. Moreover, if there is just a random noise in the picture, it doesn't make sense to run computationally intensive procedures for image reconstruction for this particular picture. Surprisingly, the vast majority of image analysis methods, both in statistics and in engineering, skip this stage and start immediately with image reconstruction.

The crucial difference of our method is that we do not impose any shape or smoothness assumptions on the \emph{boundary} of the object. This permits the detection of nonsmooth, irregular or disconnected objects in noisy images, under very mild assumptions on the object's interior. This is especially suitable, for example, if one has to detect a highly irregular non-convex object in a noisy image. This is usually the case, for example, in the aforementioned fields of automated urban analysis, cancer detection and detection of cracks in materials. Although our detection procedure works for regular images as well, it is precisely the class of irregular images with unknown shape where our method can be very advantageous.




Many modern methods of object detection, especially the ones that are used by practitioners in medical image analysis require to perform at least a preliminary reconstruction of the image in order for an object to be detected. This usually makes such methods difficult for a rigorous analysis of performance and for error control. Our approach is free from this drawback. Even though some papers work with a similar setup (see \cite{Arias-Castro_etal}), both our approach and our results differ substantially from this and other studies of the subject. We also do not use any wavelet-based techniques in the present paper.

We view the object detection problem as a nonparametric hypothesis testing problem within the class of discrete statistical inverse problems. We assume that the noise density is completely unknown, and that it is not necessarily smooth or even continuous. It is even possible that the noise distribution doesn't have a density.

In this paper, we propose an algorithmic solution for this nonparametric hypothesis testing problem. We prove that our algorithm has linear complexity in terms of the number of pixels on the screen, and this procedure is not only asymptotically consistent, but on top of that has accuracy that grows exponentially with the "number of pixels" in the object of detection. The algorithm has a built-in data-driven stopping rule, so there is no need in human assistance to stop the algorithm at an appropriate step.

In this paper, we assume that the original image is black-and-white and that the noisy image is grayscale. While our focusing on grayscale images could have been a serious limitation in case of image reconstruction, it essentially does not affect the scope of applications in the case of object detection. Indeed, in the vast majority of problems, an object that has to be detected either has (on the picture under analysis) a color that differs from the background colours (for example, in roads detection), or has the same colour but of a very different intensity, or at least an object has a relatively thick boundary that differs in colour from the background. Moreover, in practical applications one often has some prior information about colours of both the object of interest and of the background. When this is the case, the method of the present paper is applicable after simple rescaling of colour values.


The paper is organized as follows. Section 2 describes in details our statistical model, and gives a necessary mathematical introduction into the percolation theory. In Section 3, the new statistical test is introduced and its consistency is proved. This section contains the main statistical result of this paper. The power of the test and the probability of false detection for the case of small images are studied in Section 4. In addition, in Subsections 4.1 and 4.5 we describe how to select a critical cluster size for any given finite screen size. Our approach uses the fast Newman - Ziff algorithm and is completely automatic. Our main algorithm for object detection is presented in Section 5. An example illustrating the connection between the asymptotic case and the small sample case is given in Section 6. Appendix is devoted to the proof of two auxiliary estimates.


\section{Detection and percolation}

\subsection{Percolation theory}

The key observation to understand our approach to signal detection is the following central result from percolation theory \cite{Kes:82}: \\

\noindent Let $\GG$ be an infinite graph consisting of {\em sites} $s\in\GG$ and {\em bonds} between sites. The bonds determine the topology of the graph in the following sense: We say that two sites $s,s'\in\GG$ are {\em neighbors} if there is a bond connecting them. We say that a subset $C\subset\GG$ of sites is {\em connected} if for any two sites $s,s'\in C$ there are sites $s_1,...,s_n$ such that $s$ and $s_1$, $s_n$ and $s'$, and $s_k$ and $s_{k+1}$ are neighbors for all $k=1,...,n-1$. Considering {\em site percolation}\index{site percolation} on the graph $\GG$ means that we consider random configurations $\omega \in \lbrace 0,1\rbrace^{\GG}$ where the probabilites are {\em Bernoulli}
$$
\begin{array}{ll} P(\omega(s) = 1) = p, & P(\omega (s) = 0) = 1-p
\end{array}
$$
independently for each $s\in\GG$ where $0\leq p \leq 1$ is a fixed probability. If $\omega (s) =1$, we say that the site $s$ is {\em occupied}.\\

\noindent Then, under mild assumptions on the graph, there is a {\em phase transition} in the qualitative behaviour of cluster sizes. To be precise, there is a {\em critical percolation probability}\index{critical probability} $p_c$ such that for $p<p_c$ there is no infinite connected cluster and for $p>p_c$ there is one.\\

\noindent This statement and the very definition $p_c$ being the location of this phase transition are only valid for infinite graphs. We can not even speak of an infinite connected cluster for finite graphs. However, a qualitative difference of sizes of connected clusters of occupied sites can already be seen for finite graphs, say with $\vert\GG\vert = N$ sites. In a sense that will be made precise below, the sizes of connected clusters are typically of order $\log N$ for small $p$ and of order $N$ for values of $p$ close to one. This will yield a criterion to infer whether $p$ is close to zero or close to one from observed site configurations. Intuitively, for large enough values of $N$ the distinction between the two regimes is quite sharp and located very near to the critical percolation probability of an associated infinite lattice.

\subsection{The general detection problem}

Even though this paper deals mostly with a discussion of the problem for triangular lattices, we first want to sketch the problem in full generality, in order to emphasize that also other choices of the underlying lattice are possible in our processing of discretized pictures, and to make more transparent why we decided to work with six-neighborhoods in the present paper.\\


\noindent Let thus $\GG$ denote a planar graph. We think of the sites $s\in \GG$ as the pixels of a discretized image and of the graph topology as indicating neighboring pixels. We consider noisy signals of the form
\begin{equation}\label{model}
Y(s) = \bf{1}_{\GG_0} (s) + \sigma\epsilon (s)
\end{equation}
where $\bf{1}_{\GG_0}$ denotes the indicator function of a subset $\GG_0\subseteq \GG$, the noise is given by independent, identically distributed random variables $\lbrace \epsilon (s),s\in \GG\rbrace$ with $E\epsilon = 0$ and $V\epsilon = 1$, and $\sigma > 0$ is the {\em noise variance}. Thus, $\sigma^{-1}$ is a measure for the {\em signal to noise ratio}.

\begin{Definition}\label{detprob} (The detection problem)\index{detection problem} For signals of the form (\ref{model}), we consider the {\em detection problem} meaning that we construct a test for the following hypothesis
and alternative:
\begin{itemize}
\item[{$\mathbf H_0$}: ] $\GG_0=\emptyset$, i.e. there is no signal.
\item[{$\mathbf H_1$}: ] $\GG_0\neq\emptyset$, i.e. there is a signal.
\end{itemize}
\end{Definition}

\noindent{\bf Remark.} {\rm (i)} We usually think of these signals as being {\em weak} in the sense that the signal-to-noise ratio is small, i.e. $\sigma^{-1} << 1$. {\rm (ii)} Later in this paper we will specialize on the consideration of {\em symmetric} noise.\\

\noindent Our approach to the detection problem consists of translating the statements of hypothesis and alternative to statements from percolation theory. First of all, we choose a threshold $\tau\in\R$ and produce from $Y$ a {\em thresholded signal}
\begin{equation}\label{threshold}
Y_{\tau}(s) = \left\lbrace \begin{array}{ll} 1 & ,Y(s) > \tau \\ 0 & ,Y(s) \leq \tau \end{array}\right. .
\end{equation}
To adapt the terminology to percolation theory, we call a site $s\in \GG$ {\em occupied}, if $Y_{\tau}(s)=1$. Under $H_0$, the probability that a site is occupied is given by
$$
q = P(Y_{\tau}(s) = 1) = P\left(\epsilon (s) > \frac{\tau}{\sigma}\right)
$$
independently for all sites $s\in \GG$. That means, under $H_0$, the thresholded signal can be equivalently described by the occupied sites in a configuration of a {\em Bernoulli site percolation} on $\GG$ with percolation probability $q$.\\

\noindent Under $H_1$, the occupation probabilities are different for the support $\GG_0$ of the signal and its complement $\GG_1 := \GG - \GG_0$. Namely,
\begin{equation}\label{probabilities}
P(Y_{\tau} (s) = 1) = \left\lbrace\begin{array}{ll} p_0 := P\left(\epsilon > \frac{\tau - 1}{\sigma}\right) & ,s\in \GG_0 \\ p_1 := P\left(\epsilon > \frac{\tau}{\sigma} \right) & ,s\in \GG_1 \end{array}\right. .
\end{equation}
The basic idea is now to choose a threshold $\tau$ such that
\begin{equation}\label{supersub}
p_0  >  p_c  >  p_1
\end{equation}
where $p_c$ is a suitably chosen {\em critical probability} for the lattice $\GG$. That means the occupation probability is supercritical on the support of the signal and subcritical outside. \\

\begin{Remark} Note that $p_0$ and $p_1$ in the present paper denote exactly the opposite to their meaning in \cite{langovoy_report_2009-035}. We hope this causes no inconvenience for the reader, since this does not affect the main statements of this paper.
\end{Remark}

\begin{Remark} According to the intuitive picture described above, we believe that for {\em reasonably large} graphs, the critical probability $p_c$ can be chosen close or equal to the critical percolation probability of a corresponding infinite graph. See for this also the example of the triangular lattice at the end of the subsequent section.
\end{Remark}

\noindent Under mild conditions on the noise distribution which will be specified below, this can always be arranged. This is the basic observation and the starting point for our approach to the problem stated above. The idea is now to make use of the fact that the global behavior of percolation clusters is qualitatively different depending on whether the percolation probabilities are sub- or supercritical meaning that there is a detectable difference in cluster formation of occupied sites according to whether $\GG_0 = \emptyset$, or not.

\subsection[The triangular lattice and $p_c = 1/2$]{The triangular lattice and the significance of $p_c = 1/2$}\label{2_2}

From the general description of the problem formulated, some immediate questions arise:
\begin{enumerate}
\item How do we have to choose the threshold ?
\item How does the test performance depend on the choice of the underlying lattice structure for the discretized picture, i.e. the choice of $\GG$ ?
\item What can we say about the behavior of clusters of non-occupied sites in and outside of $\GG_0$ ?
\end{enumerate}
In particular, the first question is crucial. If we do not know how to properly choose the threshold to achieve super- and subcritical occupation probabilities in- and outside $\GG_0$, we have to consider scales of thresholds to determine a proper one. This would at least slow down the detection. However, the second question also indicates that we have some freedom to choose a suitable structure for the underlying lattice. It will turn out that just by choosing $\GG$ properly, i.e. as a {\em triangular lattice}, we obtain a universal answer to question {\rm (i)} and a quite satisfying result concerning question {\rm (iii)}, too. But let's start from the beginning.\\


\noindent First of all, we will summarize the conditions on the noise distribution which will be valid throughout the entire remainder of the paper.\\

\noindent{\bf Noise Properties.} For a given graph $\GG$, the noise is given by random variables $\lbrace \eps (s)\,:\,s\in\GG\rbrace$ such that
\begin{enumerate}
\item the variables $\eps (s)$ are independent, identically distributed with $E\eps=0$ and $V\eps = 1$,
\item the noise distribution is {\em symmetric},
\item the distribution of the noise is {\em non-degenerate}\index{non-degenerate noise} with respect to a critical probability $p_c$ meaning that if $F$ denotes the cumulative distribution function of the noise and we define
    $$
    \begin{array}{ll}
    m_c^+ = \inf\lbrace x\in\R\,:\,F(x)\geq 1 - p_c\rbrace , & m_c^- = \sup\lbrace x\in\R\,:\,F(x)\leq 1 - p_c\rbrace
    \end{array}
    $$
    then we have $m_c^+ = m_c^-$ where we denote the common value by $m$, and {\em either}
\begin{equation}\label{nondeg1}
F(m) > \lim_{h\to 0, h>0}F(m - h),
\end{equation}
{\em or}
\begin{equation}\label{nondeg2}
F'(m) > 0.
\end{equation}
\end{enumerate}

\noindent The reason to assume symmetry will become clear below. The reason to assume non-degeneracy is the following simple observation.

\begin{Lemma} Under the non-degeneracy condition, we can always find a threshold $\tau$ such that we have
$$
p_0 > p_c > p_1
$$
for the probabilities $p_0$, $p_1$ defined in (\ref{probabilities}). This holds independently of the value $\sigma > 0$ of the signal to noise level.
\end{Lemma}

\begin{Proof} By (\ref{probabilities}), we have $p_0 = 1 - F\left(\frac{\tau -1}{\sigma}\right)$ and $p_1 = 1 - F\left(\frac{\tau}{\sigma}\right)$. Choosing now $\tau = \sigma\,m + 1/2$, we obtain
$$
\frac{\tau -1}{\sigma} < m < \frac{\tau}{\sigma}.
$$
Under the conditions (\ref{nondeg1}) and (\ref{nondeg2}) above, that implies
$$
F\left(\frac{\tau}{\sigma}\right) > 1 - p_c > F\left(\frac{\tau -1}{\sigma}\right).
$$
\end{Proof}

\noindent Secondly, we change the terminology one last time. In the sequel, we call occupied sites {\em black} and non-occupied sites {\em white} for obvious reasons. The probability that a pixel (site) is white (not occupied) is given by
\begin{equation*}
P(Y_{\tau} (s) = 0) = \left\lbrace\begin{array}{ll} 1 - p_0 = P\left(\epsilon \leq \frac{\tau - 1}{\sigma}\right) & ,s\in \GG_0 \\ 1 - p_1 = P\left(\epsilon \leq \frac{\tau}{\sigma} \right) & ,s\in \GG_1 \end{array}\right. .
\end{equation*}
To address question {\rm (iii)} above, it would be favorable that the white pixels enjoy the same property as the black ones, just the other way round, namely that their probabilities are supercritical outside the support of the signal and subcritical inside. This will be the case if we can choose $\tau$ such that $p_0 > \max \lbrace p_c,1-p_c \rbrace \geq \min\lbrace p_c, 1- p_c\rbrace > p_1$. A situation where this can be done and where, furthermore, we get a value for $\tau$ {\em for free} is provided by the following simple but crucial observation.

\begin{Proposition}\label{triangular} Let $\GG$ be a lattice such that $p_c = 1/2$ and choose $\tau = 1/2$. If the distribution of the noise is {\em symmetric}, we have for the threshold signal $Y_\tau$:
\begin{enumerate}
\item The probability that a given pixel $s$ is black is subcritical for $s\in \GG_1$ and supercritical for $s\in \GG_0$.
\item The probability that a given pixel $s$ is white is subcritical for $s\in \GG_0$ and supercritical for $s\in \GG_1$.
\end{enumerate}
\end{Proposition}

\begin{Proof} Note first that due to the symmetry of the noise distribution, the non-degeneracy conditions (\ref{nondeg1}) and (\ref{nondeg2}) reduce to
$$\begin{array}{ll} F'(0) > 0, & \lim_{h\to 0, h >0}F(-h)<F(0).\end{array}$$
Thus, if $p_c=1/2$ and the noise is symmetric we have
$$
p_0 = P\left(\epsilon > \frac{\tau -1}{\sigma}\right) = P(\epsilon > -1/2\sigma) = P(\epsilon < 1/2\sigma) > P(\epsilon \leq 0) \geq 1/2 = p_c.
$$
Hence we also have $p_1 = P\left(\epsilon > \frac{\tau}{\sigma}\right) = P(\epsilon > 1/2\sigma) \leq 1 - P(\epsilon < 1/2\sigma) = 1 - p_0 < 1- p_c = 1/2$ and thus
$$
p_0 > \max \lbrace p_c,1-p_c \rbrace = p_c =  \min\lbrace p_c, 1- p_c\rbrace > p_1.$$
That proves the statement.
\end{Proof}

\begin{Remark} Note that, so far, the considerations are completely non-parametric with respect to the noise. Apart from the condition for expectation and variance, we do not assume anything else about the cumulative distribution function $F$ of the noise except the non-degeneracy condition and symmetry.
\end{Remark}

\noindent We round up the discussion in this section by the statement that there actually is an example for a lattice with $p_c = 1/2$.

\begin{Definition} The {\em (infinitely elongated) planar triangular lattice}\index{triangular lattice} is the infinite graph $\TT\subset\C$ in the complex plane with {\em edges} ({\em sites}) given by the elements of the additive subgroup $(S,+)\subset (\C,+)$ generated by the {\em sixth roots of unity}
$$
U := \lbrace \rho\in \C\,:\, \rho^6 = 1\rbrace.
$$
Two sites $s,s'\in S$ are connected by a {\em bond}, thus they are {\em neighbors}, if and only if their Euclidean distance in the plane is $d(s,s') =1$.
\end{Definition}

\begin{Proposition}\label{triangular2} The critical percolation probability for the {\em planar triangular lattice} is given by $p_c = 1/2$.
\end{Proposition}

\begin{Proof} See \cite{Kes:82}, p. 52 f.
\end{Proof}

\begin{Remark} Please note that this theorem holds exclusively for the infinite triangular lattice but that we rely on the assumption that for finite triangular lattices of a {\em reasonable} size, the regime of logarithmic and linear cluster sizes will be separated quite clearly and the region of probabilities $0 < p < 1$ where the regimes change will be quite sharply located around the {\em critical percolation probability} $p_c = 1/2$ of the infinite triangular lattice.
\end{Remark}

\noindent One particularly remarkable fact is that the basic receptor units ({\em ommatidia}) of an insect's eye are located at the sites of a triangular lattice. This may be caused by the requirement to have a densest possible packing of the units. However, once this discretization for visual perception is chosen, one may ask whether the analysis of the signals thus obtained use in some way or another the properties of the triangular lattice described above.

\section{The Maximum Cluster Test and consistency}\label{maxclustertest}

In this section, we will construct a test for the {\em detection problem} given in Definition \ref{detprob} above and  compute explicit upper bounds for the {\em type I and type II errors} under some mild condition on the shape of $\GG_0$, called the {\em bulk condition}. We will not focus on completeness of the proofs for which we refer to \cite{langovoy_report_2009-035} but we attempt to make transparent the basic idea of the proof: We use known statements for the infinite triangular lattice $\TT$ and transfer them to a finite lattice by arguments using a certain kind of monotonicity. The error bounds tend to zero as the lattice size aprroaches infinity. That actually provides us with a {\em consistency result}: If the image can be recorded with an unboundedly increasing resolution, the test will almost surely produce the right decision. The precise statement is given in Theorem \ref{consistency} below which is the main result of this section.


\begin{Remark}
Please note that consistency here only means that we asymptotically make the right decision. The support of the signal is not necessarily consistently reconstructed by the largest cluster. This is especially apparent when the object of interest consists of several connected components.
\end{Remark}

\noindent The setup is as follows: $\TT^{(N)}\subset \TT$ denotes the finite triangular lattice consisting of the $N^2$ sites $s\in\TT$ and bonds which are contained in the subset
$$
\lbrace z\in\C\,:\, \Re(z)\leq N + \frac{1}{2}, \Im(z) \leq \frac{\sqrt{3}}{2}N\rbrace .
$$
By consistency we mean that the test will deliver the correct decision, if the signal can be detected with an arbitrarily high resolution. To be precise, we think of the signal as a subset $G_0\subset\lbrack 0,1\rbrack^2$ and write
$$
G_0^{(N)}:=\lbrace  (N+1/2)x + iN\sqrt{3}y/2\,:\,(x,y)\in G_0\rbrace\subset \C.
$$
The model from equation (\ref{model}) is now depending on $N$, and given by
\begin{equation}\label{model_n}
Y^{(N)}(s) = \mathbf{1}_{\GG_0^{(N)}}(s) + \sigma\,\eps (s)
\end{equation}
where the sites of the subgraph are given by $\GG_0^{(N)} = \lbrace s\in\TT\,:\,s\in G_0^{(N)}\rbrace$
and the bonds of the subgraph are all bonds in $\TT$ that connect two points in $\GG_0^{(N)}$. \\

\noindent We apply now the threshold as described above, i.e. we let $\tau = 1/2$ and
$$
Y^{(N)}_{\tau}(s) = \left\lbrace\begin{array}{ll} 1 &,Y^{(N)}(s)> 1/2\\
0 & , Y^{(N)}(s)\leq 1/2 \end{array}\right. .
$$


\noindent We consider the following collection of black pixels
\begin{equation}\label{estgnull}
\hat{\GG}_0^{(N)} := \lbrace s\in\TT^{(N)}\,:\,Y^{(N)}_{\tau}(s) =1\rbrace .
\end{equation}

\noindent As a side remark, note that one can view $\hat{\GG}_0^{(N)}$ as an (inconsistent) pre-estimator of $\GG_0^{(N)}$. Now recall that we want to construct a test on the basis of this estimator for the hypotheses $\mathbf{H}_0^{(N)}: \GG^{(N)}_0 =\emptyset$ against the alternative $\mathbf{H}_1^{(N)}: \GG^{(N)}_0 \neq\emptyset$.

\begin{Definition}\label{maxclustest} {\bf (The Maximum-Cluster Test)}\index{maximum cluster test} Let $\phi (N)$ be a suitably chosen threshold depending on $N$.  Let the test statistic $T$ be the size of the largest connected black cluster $C\subset \hat{\GG}_0^{(N)}$.  We reject $\mathbf{H}_0^{(N)}$ if and only if $T \geq\phi(N)$.
\end{Definition}

\noindent For this test, we have the following consistency result under the assumption that the support of the indicator function satisfies the following very weak type of a shape constraint.


\begin{Definition}{\bf (The Bulk Condition)}\index{bulk condition} We say that the support $\GG_0^{(N)}$ of the signal {\em contains a square of side length} $\rho(N)\leq N$ if there is a site $s\in \GG_0^{(N)}$ such that $s+\TT^{(\rho(N))}\subset\GG_0^{(N)}$.
\end{Definition}

\begin{Remark} If the subset $G_0\subset \lbrack 0,1\rbrack^2$ contains a square of side length $a > 0$, the respective support $\GG_0^{(N)}$ contains squares of side lengths approximately $\rho(N)\approx a N$.
\end{Remark}

\begin{Theorem}\label{consistency} For the maximum cluster test, we have
\begin{enumerate}
\item There is some constant $K_0 > 0$ such that for $\phi (N) = K_0\log N$, we have for the {\em type I error}
$$
\lim_{N\to\infty} \alpha (N) = 0.
$$
\item Let $\phi (N)$ be as above. Let the support $\GG_0^{(N)}$ of the signal contain squares of side length $\rho (N)$. If $\rho (N) \geq K_0\log N$, we have for the {\em type II error}


$$
\lim_{N\to\infty} \beta (N) = 0.
$$
\end{enumerate}
In particular, in the limit of arbitrary large precision of sampling, the test will always produce the right detection result.
\end{Theorem}

\begin{Remark} Please note that the estimator $\hat{\GG}_0$ itself is not a consistent estimator for the support of the signal. In the best case, it is an estimator for its largest component. But we will not pursue reconstruction issues at the moment.
\end{Remark}

\noindent By the previous result, we can detect a signal correctly for virtually every noise level if we can choose an arbitrary high resolution. To prove this result, we will have to collect some facts from percolation theory to make precise the intuitive idea explained in the preceding section: {\em cluster sizes are significantly larger in the supercritical regime}. \\

\noindent We begin with the classical Aizenman-Newman theorem and transfer it step by step to finite lattices.\\

\begin{Proposition}\label{AizNew}{\bf (Aizenman-Newman Theorem)}\index{Aizenman - Newman theorem} Consider percolation with subcritical probability $p < p_c = 1/2$ on the infinite triangular lattice $\TT$. Then there is a constant $\lambda (p) > 0$ depending on $p$
such that
\begin{equation}\label{subcrit}
 P(\vert C\vert \geq n) \leq e^{-n\,\lambda(p)}
\end{equation}
for all $n\geq 1$ where $C$ denotes the black cluster containing the origin.
\end{Proposition}

\begin{Proof}See \cite{Kes:82}.\end{Proof}

\noindent Note that this result holds for the infinite triangular lattice and that we have to investigate its consequences for the finite lattices $\TT^{(N)}$. We do this by means of a classical monotonicity result known as the {\em FKG inequality}. To state this inequality, we first have to define what we mean by {\em increasing events}.

\begin{Definition}\label{increase} We can introduce a partial ordering on the set $\Omega = \lbrace 0,1\rbrace^{\TT}$ of all {\em percolation configurations} by
$$
\begin{array}{ll} \omega_1 \preceq \omega_2 & :\Longleftrightarrow \,\, \omega_1(s)\leq\omega_2(s)\,\,\mathrm{for \,\, all}\,\, s\in\TT .   \end{array}
$$
Now we say that an {\em event} $A\subset \Omega$ is {\em increasing} if we have for the corresponding indicator variable the inequality
$$
\mathbf{1}_A (\omega_1) \leq \mathbf{1}_A (\omega_2)
$$
whenever $\omega_1\preceq\omega_2$. The term {\em decreasing event} is defined analogously.
\end{Definition}

\noindent Now, the FKG inequality reads as follows.

\begin{Proposition}\label{FKG} {\bf (FKG inequality)}\index{FKG inequality} If $A$ and $B$ are both increasing (or both decreasing) events, then we have
$$
P(A\cap B) \geq P(A)\, P(B).
$$
\end{Proposition}

\begin{Proof}\cite{FoKaGi:71}\end{Proof}

\noindent We will now apply this result to compute an estimate of the amount of wrong classifications of points caused by the noise for detection of $\GG_0^{(N)}\subset \TT^{(N)}$ in a finite lattice. To be precise, we consider the event $F^{(N)}(n)$ that among the sites of a configuration $\omega\in\Omega$ that are erroneously marked black on $\TT^{(N)}$, i.e.
$$
E^{(N)}(\omega) := \lbrace s\in\TT^{(N)}-\GG_0^{(N)}\,:\,Y_{\tau}(s) = 1\rbrace
$$
and $C^{(N)} (\omega)\subset E^{(N)}(\omega)$ the largest connected cluster, then
\begin{equation}\label{FN}
F^{(N)}(n) := \lbrace \omega\in\Omega\,:\,\vert C^{(N)}\vert\geq n\rbrace .
\end{equation}
Denote now by $p_E$ the {\em error probability}
\begin{equation}\label{p_err}
p_E := P(Y_{\tau}(s) = 1\,\vert\, s\notin\GG_0^{(N)})
\end{equation}
which only depends on the noise distribution and not on $N$ or the particularly chosen site $s\notin\GG_0^{(N)}$. The consequence of the Proposition above for the finite lattice reads now as follows

\begin{Proposition}\label{expo_finite} Suppose that $0 < p_E < 1/2$ is subcritical. Let $N\geq\phi (N)$ be a threshold value with
\begin{equation}\label{quali_est}
\phi(N) = C\log N^2
\end{equation}
and $C\lambda(p_E)>1$ where $\lambda(p_E)>0$ is the value for the infinite lattice from (\ref{subcrit}). Then we have
$$
P(F^{(N)}(\phi (N)))= N^{-2(C\lambda(p_E)-1)}+ O(N^{-4(C\lambda(p_E)-1)})
$$
as $N$ tends to infinity.
\end{Proposition}

\begin{Proof} Let $s\in\TT^{(N)}$ and $C(s)$ the (possibly empty) largest connected cluster of occupied sites in $\TT$ that contains $s$. The event $\lbrace \vert C(s)\vert <n\rbrace$ is obviously {\em decreasing} for all $s\in\TT^{(N)}$. That implies by FKG inequality
\begin{eqnarray*}
P\left(\max_{s\in\TT^{(N)}} \vert C(s)\vert < n\right) &=& P\left(\cap_{s\in\TT^{(N)}} \lbrace\vert C(s)\vert < n\rbrace\right) \\
&\geq& \prod_{s\in\TT^{(N)}}P\left(\vert C(s)\vert < n\right).
\end{eqnarray*}
By (\ref{subcrit}) and translation invariance on the infinite lattice, we thus obtain
$$
P\left(\max_{s\in\TT^{(N)}} \vert C(s)\vert < n\right)
\geq \left( 1-e^{-\lambda(p_E)n}\right)^{N^2}.
$$
Hence
$$
P(F^{(N)}(n)) \leq 1 - \left( 1-e^{-\lambda(p_E)n}\right)^{N^2}
$$
Let now $n = 2C\log N = \phi (N)$ with $C\lambda (p_E) > 1$. Then, by (\ref{asymptotics}), we obtain
\begin{equation*}
P(F^{(N)}(\phi(N))) = N^{-2(C\lambda(p_E)-1)}+ O(N^{-4(C\lambda(p_E)-1)}).
\end{equation*}
\end{Proof}

\begin{Definition}\label{left_right} Let $\TT^{(N)}\subset \TT$ be as above. A subset $\pi = \lbrace s_1,...,s_n\rbrace$ of black sites $s_k\in\TT^{(N)}$ with
\begin{enumerate}
\item $s_k$ and $s_{k+1}$ are neighboring sites for all $k=1,...,n-1$,
\item $0\leq \Re (s_1)\leq 1/2$,
\item $N\leq \Re(s_n)\leq N+1/2$,
\end{enumerate}
is called a {\em left-right crossing}.
\end{Definition}

\begin{Proposition}\label{supercrit} Consider site percolation on $\TT^{(N)}\subset\TT$ with supercritical percolation probability $p > 1/2$. Denote by $A_N$ the event that there is some left-right crossing in $\TT^{(N)}$. Then there is a constant $D(p) > 0$ such that
$$
P(A_N) \geq 1- Ne^{-D(p)N}.
$$
\end{Proposition}

\begin{Proof} The proof is analogous to the one of the corresponding statement for bond percolation in \cite{Gri:99}. See \cite{langovoy_report_2009-035} for more information.
\end{Proof}

\begin{Proof}{\bf (of Theorem \ref{consistency})} {\rm(i)} By Proposition \ref{triangular}, the probability $p=p_E$ for a pixel to be black is subcritical under the null hypothesis. Let $\lambda(p)$ be the exponential factor for the infinite lattice in equation (\ref{subcrit}). Choosing now $K_0 = 2C$ with $C\lambda(p_E)>1$ means that by Proposition \ref{expo_finite} the type I error probability tends to zero as $N$ tends to infinity.  \par

\noindent {\rm (ii)} Let the alternative be true and  $\GG_0^{(N)}\neq \emptyset$ contain a square of side length $\rho (N)$. Again by Proposition \ref{triangular}, the probability $p_B>1/2$ that a site is (correctly) marked black is supercritical for sites inside the square. Therefore, the probability to find a left-right crossing $\pi$ is given by Proposition \ref{supercrit} to be
$$
P(A_{\rho (N)}) \geq 1- \rho(N)e^{-D(p_B)\rho(N)}.
$$
But every left-right crossing is a connected cluster of size at least $\rho(N)$. Hence
$$
P(T \geq\rho(N)) \geq 1- \rho(N)e^{-D(p_B)\rho(N)},
$$

\noindent and together with $\rho(N)\geq K_0\log N$ this implies that
$$
\lim_{N\to\infty}P(T \geq K_0\log N) \geq \lim_{N\to\infty}P(T \geq \rho(N)) = 1
$$
and thus $\lim_{N\to\infty}\beta (N) = 0$.
\end{Proof}

\noindent Finally, as an addition to Proposition \ref{expo_finite}, we prove that under the same conditions, a slightly changed arguments yields an exponential estimate for the tail probabilities of the cluster size distribution.

\begin{Proposition}\label{good_exponential_tails} Suppose that $0 < p_E < 1/2$ is subcritical. Let $N\geq\phi (N)$ be a threshold value with
\begin{equation*}
\phi(N) = C\log N^2
\end{equation*}
and $C\lambda(p_E)>1$ where $\lambda(p_E)>0$ is the value for the infinite lattice from (\ref{subcrit}). Then there exists a constant $K^{(N)}> 0$ such that
$$
P(F^{(N)}(n))\leq e^{-K^{(N)}n}
$$
for all $n\geq \phi (N)$.
\end{Proposition}

\begin{Proof} Let $s\in\TT^{(N)}$ and $C(s)$ the (possibly empty) largest connected cluster of occupied sites in $\TT$ that contains $s$. The event $\lbrace \vert C(s)\vert <n\rbrace$ is obviously {\em decreasing} for all $s\in\TT^{(N)}$. That implies by FKG inequality
\begin{eqnarray*}
P\left(\max_{s\in\TT^{(N)}} \vert C(s)\vert < n\right) &=& P\left(\cap_{s\in\TT^{(N)}} \lbrace\vert C(s)\vert < n\rbrace\right) \\
&\geq& \prod_{s\in\TT^{(N)}}P\left(\vert C(s)\vert < n\right).
\end{eqnarray*}
By (\ref{subcrit}) and translation invariance on the infinite lattice, we thus obtain
$$
P\left(\max_{s\in\TT^{(N)}} \vert C(s)\vert < n\right)
\geq \left( 1-e^{-\lambda(p_E)n}\right)^{N^2}.
$$
Hence
$$
P(F^{(N)}(n)) \leq 1 - \left( 1-e^{-\lambda(p_E)n}\right)^{N^2}
$$
Let now $n\geq 2C\log N$ with $C\lambda (p_E) > 1$. Then, by essentially the same calculation as for (\ref{estimate}), we obtain
\begin{eqnarray*}
P(F^{(N)}(n)) &\leq& N^2 e^{-\lambda(p_E)n} \left(1 + e^{-\lambda(p_E)n}\right)^{N^2-1}\\
&\leq& N^2  \left(1 + N^{-2C\lambda(p_E)}\right)^{N^2-1} e^{-\lambda(p_E)n}.
\end{eqnarray*}
Now, $n\geq C\log N^2$ and thus $\log N^2 \leq n/C$. That implies by $C\lambda > 1$ and
$$
\lim_{N\to\infty} \left(1 + N^{-2C\lambda(p_E)}\right)^{N^2-1} =1
$$
that for $N$ large enough, we have
$$
\left(1 + N^{-2C\lambda(p_E)}\right)^{N^2-1} \leq e^{\kappa n}
$$
and therefore
$$
\log N^2 +\kappa - \lambda (p_E)n \leq \left(\frac{1}{C} - \lambda (p_E) \right) n = -K^{(N)}n
$$
with $K^{(N)} := \lambda(p_E) - \kappa-C^{-1} > 0$.
\end{Proof}

This proposition helps to strengthen Theorem \ref{consistency} and to derive the actual rates of convergence for both types of testing errors. It is a remarkable fact that both types of errors in our method tend to zero exponentially fast in terms of the size of the object of interest.

\begin{Theorem}\label{rates}
Suppose assumptions of Theorem \ref{consistency} are satisfied. Then there are constants $C_1 > 0$, $C_2 > 0$ such that\smallskip

\begin{enumerate}

  \item The type I error of the maximum cluster test does not exceed 
  \[
  \alpha (N) \leq \exp(-C_2 \phi (N))
  \] 
  
  \noindent for all $N > \phi (N)$.\medskip

  \item The type II error of the maximum cluster test does not exceed 
  \[
  \beta (N) \leq \exp(-C_1 \rho (N))) \,.
  \]
  
   \noindent for all $N > \rho (N)$.
\end{enumerate}

\end{Theorem}

\begin{Proof}
Analogously to the proof of Theorem \ref{consistency}, only replacing Proposition \ref{expo_finite} by Proposition \ref{good_exponential_tails} and using the inequality from Proposition \ref{supercrit} in a slightly modified fashion, just as was done in the proof of Theorem 1 in \cite{langovoy_report_2009-035}.
\end{Proof}

\noindent In other words, the maximum cluster test has power that goes to one exponentially, and the false detection rate that goes to zero exponentially, comparatively with the size of the object of interest.

\section{Power, significance level and critical size of clusters}

\subsection{Type I error and the critical size of clusters}

By the consistency result obtained above, we see that the maximum cluster test leads almost surely to the right test decision if we can sample the incoming signal with arbitrary precision. However, as frequently, the asymptotic result will not determine the exact threshold $\phi (N)$ for a given {\em finite} value of $N$ when we assume some level of significance $\alpha > 0$. Since the significance level is based on the determination of the type one error, we can use the fact the under the null hypothesis $\mathbf{H}_0 : f = 0$, the black clusters on the graph are distributed according to a plain percolation with probability of a site to be black equal to $p_E$. In that case, we can use the {\em Newman - Ziff algorithm} which is described in \cite{Wittich_Langovoy_R} to effectively simulate the distribution of the size of the largest cluster. \\

\noindent In the sequel, we present the results of these simulations for the corresponding critical regions of the maximum cluster size statistic $T$ for a triangular lattice with $55\times 55 = 3025$ sites. We consider {\em significance levels} of $\alpha = 0.05$ and $\alpha = 0.01$ and different possible values for $p_E$. We obtain the following tables for the lower bounds of the critical regions $\lbrace T \geq c_0\rbrace$ depending on values for $p_E$.

\vspace{0.1cm}

\begin{center}
\begin{tabular}{|l|c|c|c|c|c|c|c|c|c|}\hline $\mathrm{p_E}$ & 0.1 & 0.2 & 0.3 & 0.4 & 0.42 & 0.44 & 0.46 & 0.48 & 0.5\\\hline
$\mathrm{c_0}$ & 7  & 19 &  49 & 186 & 262 & 395 & 597 & 891 & 1184 \\\hline
\end{tabular}

\vspace{0.1cm}

{\em Boundary of critical regions for $\alpha = 0.05$}
\end{center}

\vspace{0.1cm}

\begin{center}
\begin{tabular}{|l|c|c|c|c|c|c|c|c|c|}\hline $\mathrm{p_E}$ & 0.1 & 0.2 & 0.3 & 0.4 & 0.42 & 0.44 & 0.46 & 0.48 & 0.5\\\hline
$\mathrm{c_0}$ & 9 &  23  & 62 & 247 & 352 & 524 & 765 & 1058 & 1302 \\\hline
\end{tabular}

\vspace{0.1cm}

{\em Boundary of critical regions for $\alpha = 0.01$}
\end{center}

\vspace{0.1cm}

\noindent These values are obtained from the quantiles of the simulated cluster size distributions according to the description given in \cite{Wittich_Langovoy_R}. The site probabilities for which the simulations are run, are $p_E$ = 0.1, 0.2,  0.3,  0.4, 0.42,  0.44,  0.46,  0.48,  0.5,  0.52,  0.54,  0.56,  0.58,  0.6,  0.7,  0.8,  and 0.9. They are displayed in the following graphics. \\

\begin{figure}[!h]\label{cdf}
\begin{center}
\fbox{
\begin{tabular}{cc}\includegraphics[width=0.40\textwidth]{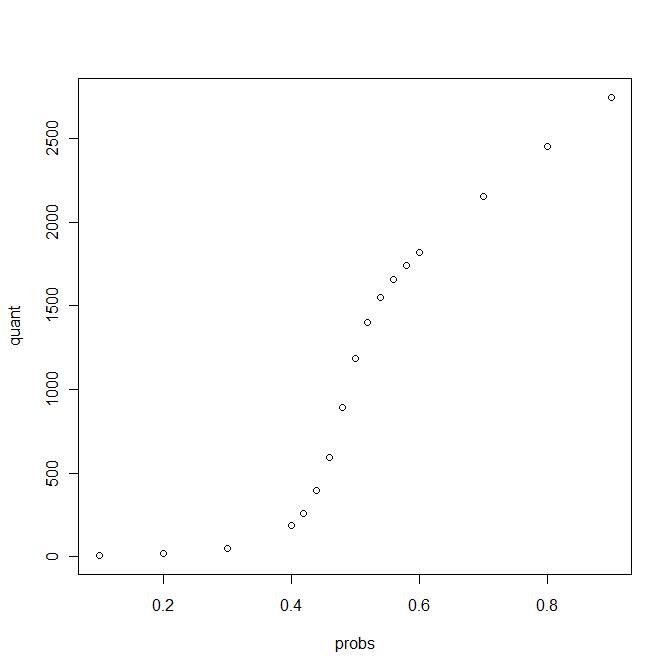} & \includegraphics[width=0.40\textwidth]{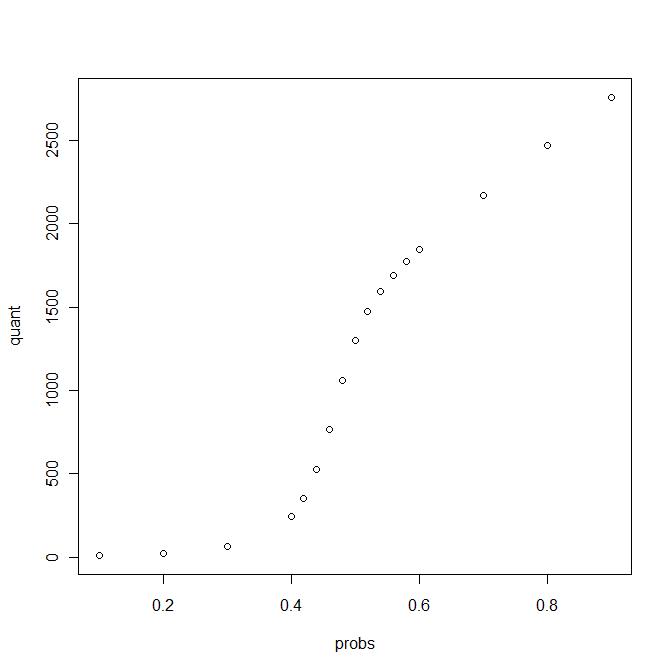}\\
\end{tabular}
}
\end{center}
\caption{0.95- and 0.99-quantiles of maximum cluster size as function of $p_E$}
\end{figure}

\subsection{Type II error -- monotonicity}

To compute upper and lower estimates of the type II probability, we will use the bulk condition together with the simple fact that the error probability of the test is monotonous in the support of the signal. To be precise, we have the following statement.

\begin{Lemma} Let $\GG_0\subset\GG_0'\subset\GG$ be two different supports for signals $Y,Y'$ according to equation (\ref{model}). Then, for any given level of significance, the type II error for signal $Y$ is larger than for $Y'$.
\end{Lemma}

\begin{Proof} For every configuration $\omega\in\lbrace 0,1\rbrace^{\GG}$ after thresholding, we have using the probabilities $p_1 > p_0$ from (\ref{supersub}) and the notation $\GG_0$, $\GG_1$ from(\ref{probabilities})
\begin{eqnarray*}
P(\omega) &=&  p_0^{\vert\omega^{-1}(1)\cap\GG_0\vert}(1 - p_0)^{\vert\omega^{-1}(0)\cap\GG_0\vert}p_1^{\vert\omega^{-1}(1)\cap\GG_1\vert}(1-p_1)^{\vert\omega^{-1}(0)\cap\GG_1\vert}\\
&=&  \frac{p_1^{\vert\omega^{-1}(1)\cap\GG_1\vert - \vert\omega^{-1}(1)\cap\GG_1'\vert}(1 - p_0)^{\vert\omega^{-1}(0)\cap\GG_0\vert - \vert\omega^{-1}(0)\cap\GG_0'\vert}}{p_0^{\vert\omega^{-1}(1)\cap\GG_0'\vert -\vert\omega^{-1}(1)\cap\GG_0\vert}(1-p_1)^{\vert\omega^{-1}(0)\cap\GG_1'\vert - \vert\omega^{-1}(0)\cap\GG_0\vert}}\times\\
&& \times  p_0^{\vert\omega^{-1}(1)\cap\GG_0'\vert}(1 - p_0)^{\vert\omega^{-1}(0)\cap\GG_0'\vert}p_1^{\vert\omega^{-1}(1)\cap\GG_1'\vert}(1-p_1)^{\vert\omega^{-1}(0)\cap\GG_1'\vert}\\
&\leq& p_0^{\vert\omega^{-1}(1)\cap\GG_0'\vert}(1 - p_0)^{\vert\omega^{-1}(0)\cap\GG_0'\vert}p_1^{\vert\omega^{-1}(1)\cap\GG_1'\vert}(1-p_1)^{\vert\omega^{-1}(0)\cap\GG_1'\vert}\\
&=:& P'(\omega)
\end{eqnarray*}
by $p_0 > p_1$ and
$$
\vert\omega^{-1}(1)\cap\GG_1\vert - \vert\omega^{-1}(1)\cap\GG_1'\vert = \vert\omega^{-1}(1)\cap\GG_0'\vert -\vert\omega^{-1}(1)\cap\GG_0\vert
$$
 and where $P$ denotes the the probability of the configuration associated to $Y$ and $P'$ the probability associated to $Y'$. That implies for the probability that the that the signal is not detected $P(T < n_0) \geq P'(T < n_0)$ where $n_0$ is a cluster size that is determined on the basis of the null hypothesis and the level of significance alone and does therefore not depend on the signal. That implies the statement.
\end{Proof}

Thus, under the {\em bulk condition}, we obtain a conservative upper estimate for the type II error, if we simulate the type II error for the square contained in the support of the signal and a lower estimate if we compute it for a square containing the support. In the latter case, we will consider here the best case scenario of a signal that is supported by all of $\GG$.

\subsection{Type II error -- lower bound for small images}

In this paragraph we study numerical performance of our procedure in the case of small images. For illustrative purposes, we concentrate on the case $N = 55$, which is important for applications of our method in insect vision (see \cite{Wittich_Langovoy_Insect_Vision}). To estimate the type II error in this case, we compute the probability that a constant signal of the form $f(s) = a > 0$ (which is therefore positive on the whole graph) is not detected. So we actually consider the alternative $\mathbf{H}_1: f = a > 0$. The probability that a given pixel is marked black is therefore supercritical given by $p_B > 1/2$. From the simulated distribution of the maximum cluster sizes, we thus obtain for the type II error the subsequent tables by just taking the value of the cumulative distribution function at the boundaries of the respective critical regions.

\vspace{0.2cm}

\begin{center}
\begin{tabular}{|l|c|c|c|c|c|c|c|c|c|}\hline $\begin{array}{c} \mathrm{p_E}\to\\  \mathrm{p_B}\downarrow  \end{array}$ & 0.1 & 0.2 & 0.3 & 0.4 & 0.42 & 0.44 & 0.46 & 0.48 & 0.5\\\hline
0.52 &   &  &   & $10^{-6}$ & 0.0005 & 0.01 & 0.09 & 0.31 & 0.66 \\\hline
0.54 &   &  &   & $5\, 10^{-9}$ & $10^{-5}$ & 0.0007 & 0.01 & 0.08 & 0.23 \\\hline
0.56 &   &  &   &  & $2\,10^{-9}$ & $6\,10^{-7}$ & 0.0006 & 0.01 & 0.05 \\\hline
0.58 &   &  &   &  &  &  & $7\,10^{-6}$ & 0.0009 & 0.006 \\\hline
0.6 &   &  &   &  &  &  & $8\, 10^{-9}$ & $4\,10^{-5}$ & 0.0005 \\\hline
0.7 &   &  &   &  &  &  &  &  &  \\\hline
0.8 &   &  &   &  &  &  &  &  &  \\\hline
0.9 &   &  &   &  &  &  &  &  &  \\\hline
\end{tabular}

\vspace{0.1cm}

{\em Type II error for $\alpha = 0.05$, all omitted values are less than $10^{-9}$}
\end{center}

\vspace{0.1cm}

\begin{center}
\begin{tabular}{|l|c|c|c|c|c|c|c|c|c|}\hline $\begin{array}{c} \mathrm{p_E}\to \\  \mathrm{p_B}\downarrow \end{array}$ & 0.1 & 0.2 & 0.3 & 0.4 & 0.42 & 0.44 & 0.46 & 0.48 & 0.5\\\hline
0.52 &   &  &   & 0.0003 & 0.006 & 0.05 & 0.21 & 0.49  & 0.83 \\\hline
0.54 &   &  &   & $10^{-5}$ & 0.0003 & 0.005 & 0.04 & 0.15 & 0.39 \\\hline
0.56 &   &  &   &  & $2\,10^{-7}$ & 0.0002 & 0.004 & 0.03 & 0.09 \\\hline
0.58 &   &  &   &  &  & $2\,10^{-6}$ & 0.0002 & 0.003 & 0.01 \\\hline
0.6 &   &  &   &  &  &  & $10^{-5}$ & 0.0003 & 0.0009 \\\hline
0.7 &   &  &   &  &  &  &  &  &  \\\hline
0.8 &   &  &   &  &  &  &  &  &  \\\hline
0.9 &   &  &   &  &  &  &  &  &  \\\hline
\end{tabular}

\vspace{0.1cm}

{\em Type II error for $\alpha = 0.01$, all omitted values are less than $10^{-9}$}
\end{center}

\vspace{0.2cm}

\noindent Of course, we would ideally prefer to simulate the type II error for signals containing certain squares such as used for the consistency result. To be precise, we would like to consider the case where the real signal is of the form $a\,\mathbf{1}_Q$ where $Q\subset \TT^{(N)}$ is a given square and $a > 0$. However, we can not use the Newman - Ziff algorithm directly for that since the site probabilities are inhomogeneous. To find an efficient simulation algorithm for the type II error in this case is therefore work in progress.

On the other hand, the above two tables show that our testing procedure can have very good power already for rather small images. Since it follows from Theorem \ref{rates} that both error probabilities tend to zero exponentially as the screen resolution increases, both the power and the level of the Maximum Cluster Test will improve rapidly even with a small growth of the image resolution.

\subsection{Type II error -- upper bound for small images}\label{upbound}

To simulate the type II error for a signal supported on a square sublattice, we use the modified version of the Newman - Ziff algorithm described in \cite{Wittich_Langovoy_R}. The simulations are done for a $15\times 15$ square lattice that is basically situated in the center of a $55\times 55$ lattice. We expect that the simulated probabilities do depend on the location of the support but for us the simulated type II errors in this section only serve as a proof of principle. To further explore the type II error will take extensive simulations for different shapes and locations of the support, but we will not pursue that here.

\subsection{Simulation of the maximum cluster size distribution}

All the experiments described so far concerned the maximum cluster size distribution of a triangular lattice $\TT{(55)}$ with $55\times 55 = 3025$ sites for different site probabilities of $p_E$ = 0.1, 0.2,  0.3,  0.4, 0.42,  0.44,  0.46,  0.48,  0.5,  0.52,  0.54,  0.56,  0.58,  0.6,  0.7,  0.8,  and 0.9. The results of these simulations are presented in Figures 2 and 3 below. We used the R-implementation {\tt simulation} of the Newman - Ziff algorithm described in \cite{Wittich_Langovoy_R}.

\begin{figure}
\begin{center}
\fbox{
\begin{tabular}{cc}
\includegraphics[width=0.30\textwidth]{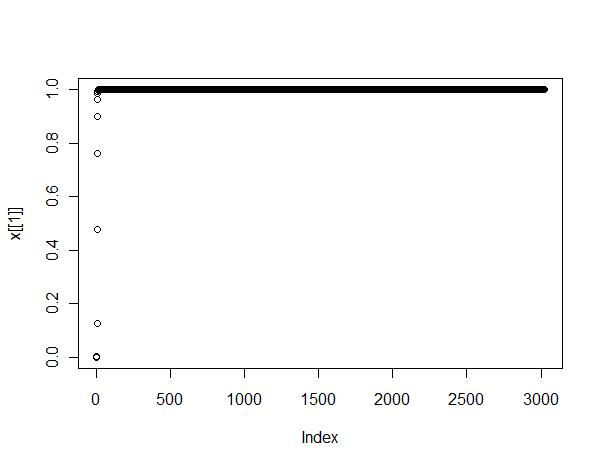} & \includegraphics[width=0.30\textwidth]{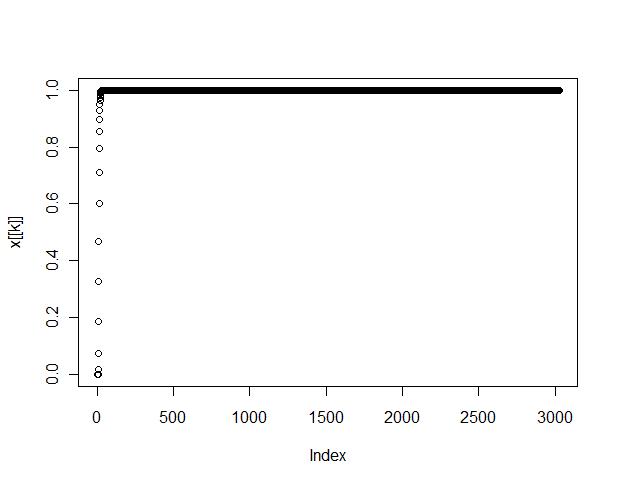}\\
\includegraphics[width=0.30\textwidth]{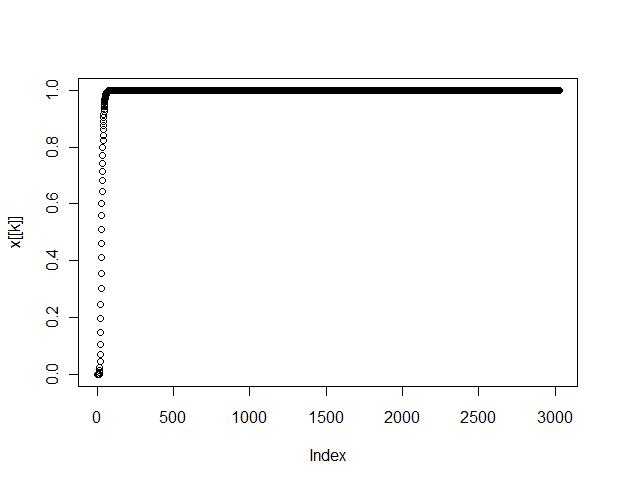} & \includegraphics[width=0.30\textwidth]{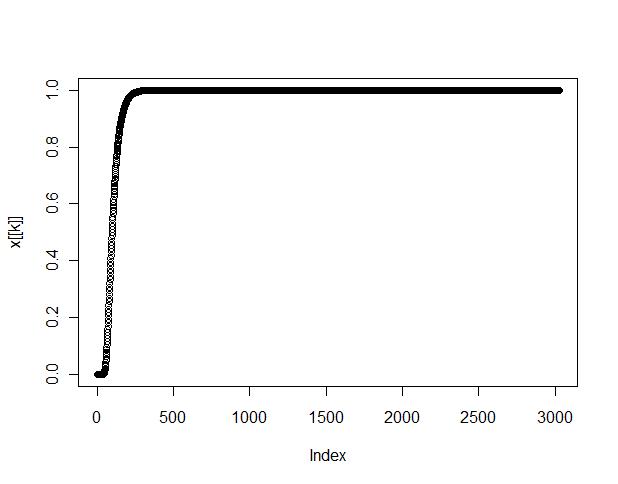}\\
\includegraphics[width=0.30\textwidth]{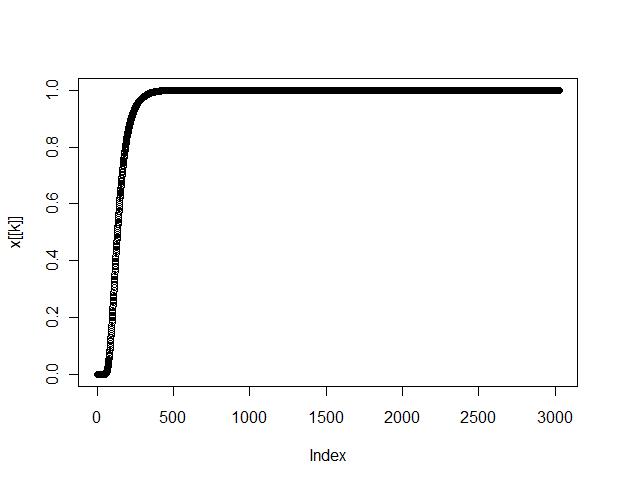} & \includegraphics[width=0.30\textwidth]{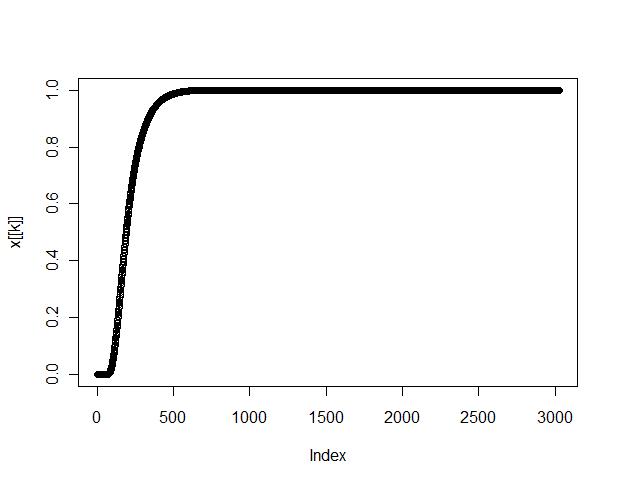}\\
\includegraphics[width=0.30\textwidth]{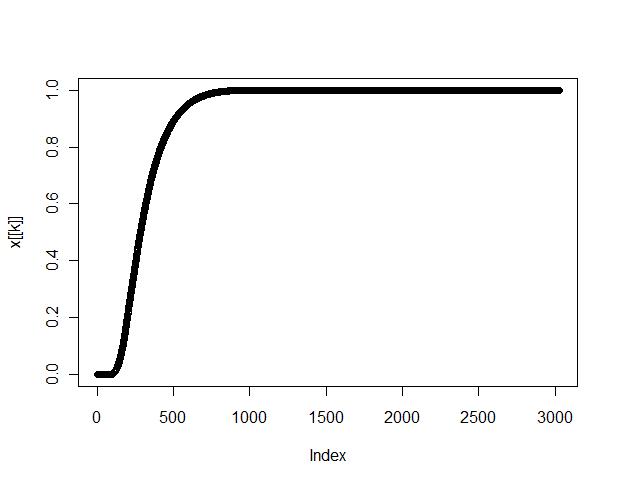} & \includegraphics[width=0.30\textwidth]{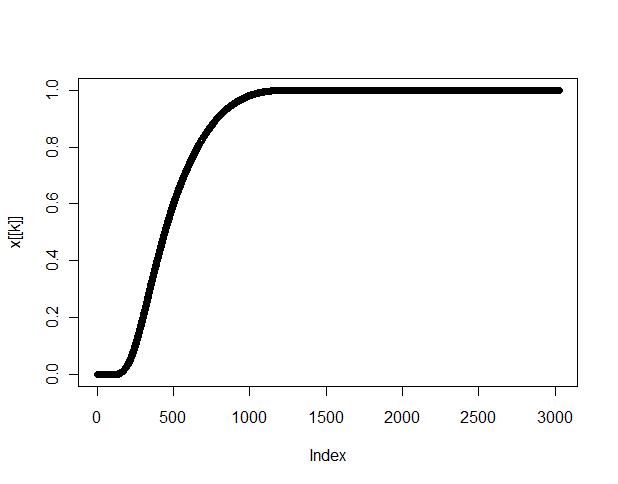}\\
\includegraphics[width=0.30\textwidth]{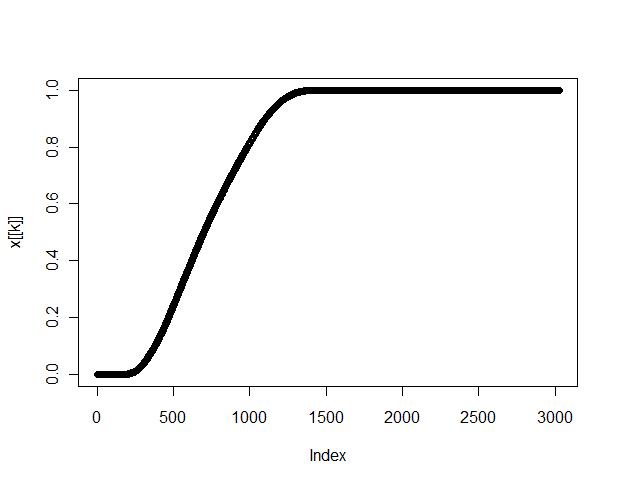} & \includegraphics[width=0.30\textwidth]{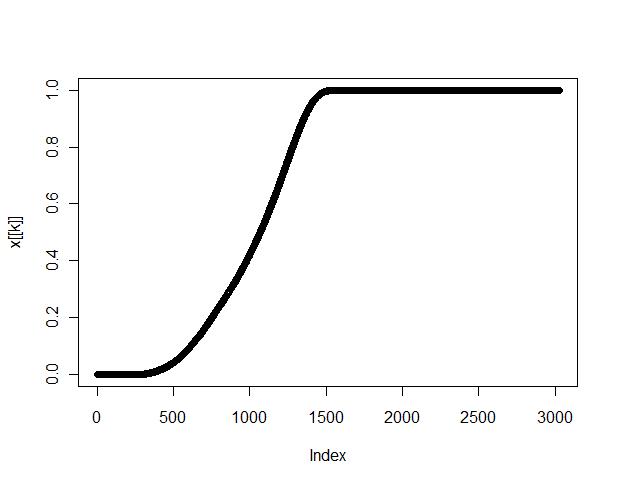}
\end{tabular}
}
\end{center}
\caption{Maximum cluster size distribution for $p_E$ = 0.1, 0.2, 0.3, 0.4, 0.42, 0.44, 0.46, 0.48, 0.5, 0.52}
\end{figure}

\begin{figure}[!t]
\begin{center}
\fbox{
\begin{tabular}{cc}
\includegraphics[width=0.30\textwidth]{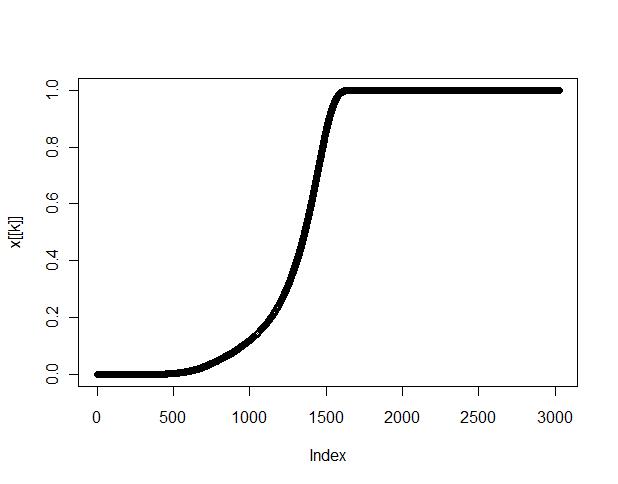} & \includegraphics[width=0.30\textwidth]{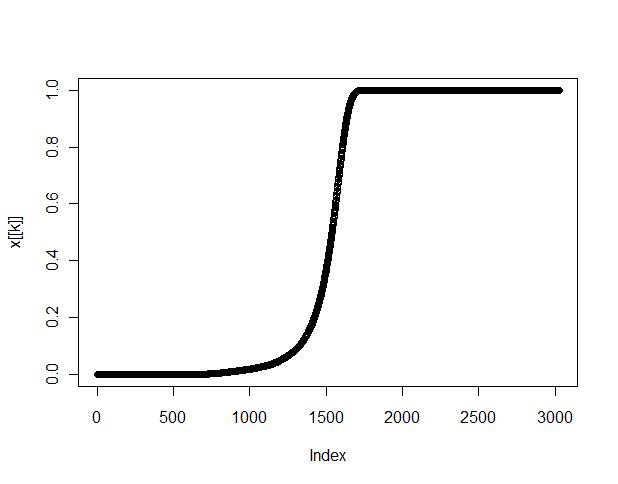}\\
\includegraphics[width=0.30\textwidth]{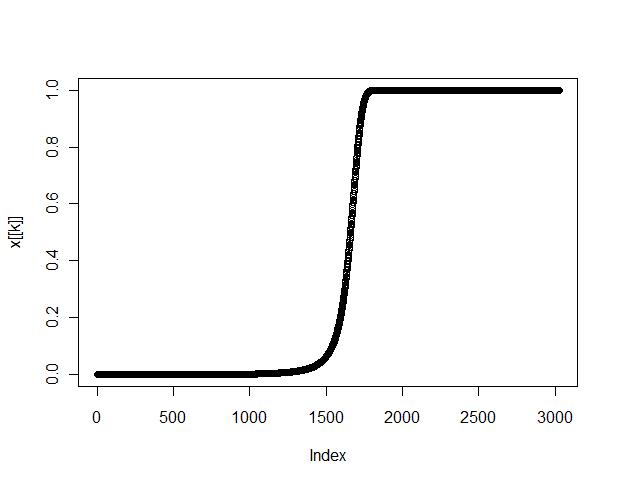} & \includegraphics[width=0.30\textwidth]{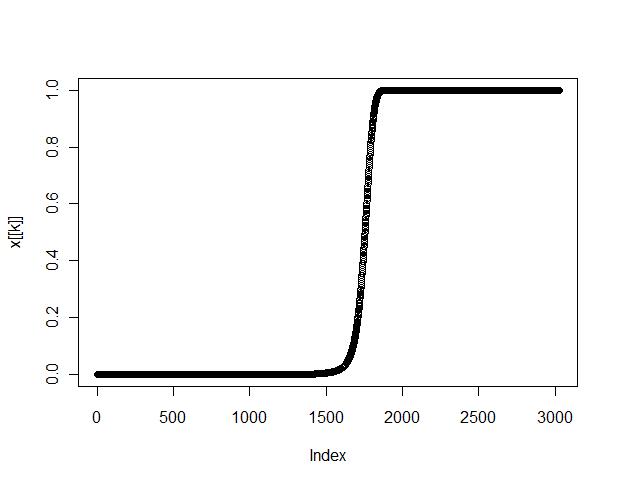}\\
\includegraphics[width=0.30\textwidth]{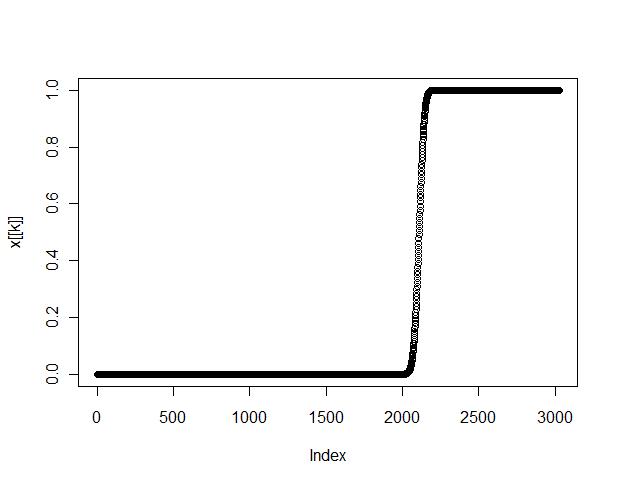} & \includegraphics[width=0.30\textwidth]{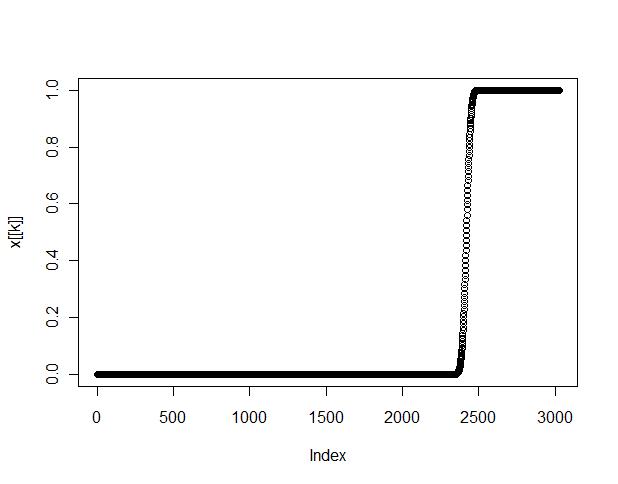}\\
\includegraphics[width=0.30\textwidth]{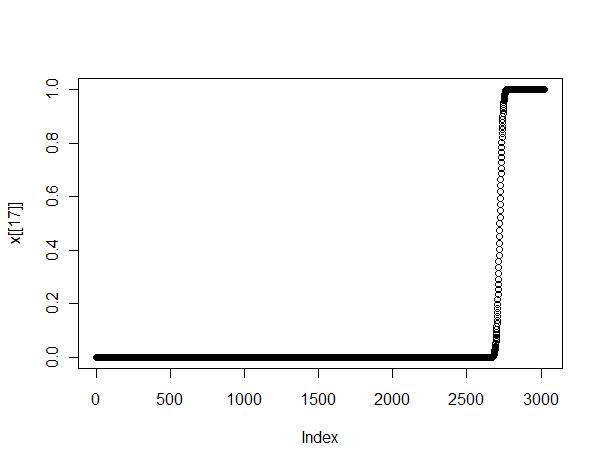} & \\
\end{tabular}
}
\end{center}
\caption{Maximum cluster size distribution for $p_E$ = 0.54, 0.56, 0.58, 0.6, 0.7, 0.8, 0.9}
\end{figure}



\section[Detection algorithm]{Detection algorithm}

Once a threshold $\tau$ is fixed, we will base the test decision (i.e. {\em we detect a signal or not}) on the size of the maximal black cluster. It is therefore certainly an important point that there is a cluster search algorithm, the {\em Depth First Search}\index{Depth First Search algorithm} algorithm from \cite{Tar:72} which is explained in details and implemented in R in \cite{Wittich_Langovoy_R}, which is quite effective. That means, it is linear in the number of pixels. This is stated below. Please note that the R-implementation is rather slow compared to an implementation in C.\\

\noindent We now describe the detection algorithm and state the complexity result. The algorithm consists of the following steps:
\begin{enumerate}
\item Perform a $\tau$-thresholding of the noisy picture $\hat{Y}$.
\item Run a {\em depth first search} algorithm on the graph of the thresholded signal until {\em either} a black cluster of size $\vert C \vert \geq c_0$ is found, {\em or} all black clusters are found.
\item If a black cluster of size $\vert C\vert \geq c_0$ was found, report that a signal was detected, otherwise do not reject $\mathbf{H}_0$.
\end{enumerate}
The complexity result follows now from the linear complexity of Tarjan's algorithm (see \cite{Tar:72}) and is given by the following statement.

\begin{Theorem} The algorithm terminates in $O(N^2)$ steps, i.e. it is linear in the number of pixels.
\end{Theorem}

\begin{Proof} See \cite{langovoy_report_2009-035}, Theorem 1. \end{Proof}


\section{Application of asymptotic results to small images}

Finally, we want to consider one (hopefully) instructive example. We assume that $\eps \sim N(0,1)$ is standard normally distributed and that the true signal is given by the constant $a > 0$. We choose for the test the threshold $\tau > 0$ and obtain
$$
p_E(\tau) = P(\eps > \tau/\sigma) = 1 - \Phi(\tau/\sigma) < 1/2.
$$
Given the true signal, the correct probability that a site is marked black is given by
$$
p_B(\tau, a) = P(\eps > (\tau - a)/\sigma) = \Phi ((a - \tau)/\sigma).
$$
The first observation is that $p_B$ is supercritical, only if $\tau\in\lbrack 0, a)$. This may help us to develop an algorithm to find good threshold values in the following way.
\begin{enumerate}
\item Start with some threshold value $\tau_1 > 0$ and perform the algorithm described above.
\item If $\mathbf{H}_0$ is not rejected, proceed with a smaller value $\tau_2$.
\item Repeat step 2 until you reject $\mathbf{H}_0$ or $\tau_n < \tau_0$ where $\tau_0$ is some minimal threshold given in advance.
\end{enumerate}
To use such a minimal threshold makes sense due to the uncertainty relation, which in the context of a finite lattice with simulated maximum cluster size distribution corresponds to the fact that the type II error can get large. To be precise, for a given (small) $a>0$ we have, because supercritical behavior can only be achieved for $\tau < a$, that
$$
1/2 > p_E (\tau) > 1 - \Phi (a/\sigma),\,\,\,1/2 < p_B (\tau) < \Phi (a/\sigma).
$$

\begin{Remark} From the first sight, it seems obvious that given a signal strength $a > 0$, the threshold $\tau = a/2$ is optimal in separating the two probabilities $p_E$ and $p_B$. But due to the asymmetry of the distribution functions around the critical probability $1/2$ visible in Figure \ref{cdf}, this may not be the case. Moreover, this is only the case for certain symmetric lattices and symmetric noise distributions.
\end{Remark}

\noindent However, the inequalities above already yield that
$$
1 - \Phi (a/\sigma) < p_E (\tau) < 1/2 < p_B(\tau) < \Phi (a/\sigma)
$$
and that means for instance if $\alpha = 0.01$ and the signal to noise ratio is $\rho = a/\sigma = 0.05$, then $p_B \approx 0.52$, $p_E\approx 0.48$, we may infer from the second table about type II errors that $\beta \approx 0.49$. It may therefore not really make sense to consider threshold values below $\tau _0 = 0.05\times \sigma$. 


\smallskip
\noindent {\bf Acknowledgments.} The authors would like to thank Laurie Davies and Remco van der Hofstad for helpful discussions. \\

\bibliographystyle{plainnat}

\bibliography{papiere}

\bigskip
\noindent {\bf Appendix.}\\



We will shortly explain how to obtain an asymptotic expression of the binomial sum in the proof of Proposition \ref{expo_finite}. 

\begin{Lemma} For $\vert x\vert < 1$, $1\leq n\in\N$, we have
\begin{equation}\label{estimate}
\vert(1+x)^n-1-nx\vert\leq\frac{n(n-1)}{2}\vert x\vert^2 (1+\vert x \vert)^{n-2}.
\end{equation}
\end{Lemma}

\begin{Proof} We have
\begin{eqnarray*}
\vert(1+x)^n-1-nx\vert &=& \sum_{k=2}^n\left(\begin{array}{c}n \\ k\end{array}\right)\vert x\vert^k\\
&=& \frac{n(n-1)}{2}\vert x\vert^2\sum_{k=2}^n\frac{\left(\begin{array}{c}n \\ k\end{array}\right)}{\left(\begin{array}{c}n \\ 2\end{array}\right)}\vert x\vert^{k-2}\\
&\leq& \frac{n(n-1)}{2}\vert x\vert^2\sum_{k=0}^{n-2}\left(\begin{array}{c}n-2 \\ k\end{array}\right)\vert x\vert^{k}\\
&=&\frac{n(n-1)}{2}\vert x\vert^2 (1+\vert x \vert)^{n-2}.
\end{eqnarray*}
\end{Proof}

\begin{Lemma} Let $z_N = e^{-\lambda f(N)}$ and $f(N)=C\log (N^2) = 2C\log(N)$ where $\lambda > 0$ and $\lambda C > 1$. Then, as $N$ tends to infinity, we have
\begin{equation}\label{asymptotics}
1 - (1-z_N)^{N^2} = N^{2(1-C\lambda) }+ O\left( N^{4(1-C\lambda)}\right).
\end{equation}
\end{Lemma}

\begin{Proof} By inequality (\ref{estimate}) above
$$
\vert(1-z_N)^{N^2} - 1 + N^2 z_N\vert\leq \frac{N^2(N^2-1)}{2} z_N^2 (1+\vert z_N\vert)^{N^2-2}
$$
By $C\lambda > 1$, we obtain
$$
\lim_{N\to\infty} (1+\vert z_N \vert)^{N^2-2} = \lim_{N\to\infty} \left(1+\frac{1}{N^{2C\lambda}}\right)^{N^2} = 1.
$$
Thus, for $N>N_0$ large enough, there are constants $K^* > K>1$ such that
$$
\vert(1-z_N)^{N^2} - 1 + N^2 z_N\vert\leq K\frac{N^2(N^2-1)}{2} z_N^2 \leq K^* N^{4(1-C\lambda)}.
$$
Together with
$$
N^2 z_N = N^{2(1-C\lambda) },
$$
that proves the assertion.
\end{Proof}


\end{document}